\documentstyle[12pt]{article}
\def\e{{\varepsilon}}

\def\a{{\alpha}}

\def\lar{{\leftarrow}}

\def\be{\begin{equation}}
\def\ee{\end{equation}}

\newtheorem{thm}{\hspace{\parindent}Theorem}[section]
\newtheorem{lem}[thm]{\hspace{\parindent}Lemma}
\newtheorem{prop}[thm]{\hspace{\parindent}Proposition}
\newtheorem{dfn}[thm]{\hspace{\parindent}Definition}

\newcommand{\norm}[1]{\left\| #1 \right\|}
\newcommand{\til}[1]{\widetilde{#1}}
\newcommand{\ov}[1]{\overline{#1}}

\date{}
\author{V. M. Manuilov,\quad A. S. Mishchenko}
\title{Relations between asymptotic and Fredholm representations%
\footnote{This research was partially supported by
RFBR (grant No 96-01-00276) and by DFG.}}

\begin{document}
\maketitle

\begin{abstract}
We prove that for matrix algebras $M_n$ there exists a monomorphism
$$
\left({\textstyle \prod}_n M_n/\oplus_n M_n\right) \otimes C(S^1)
\longrightarrow {\cal Q}
$$
into
the Calkin algebra which induces an isomorphism of the $K_1$-groups. As a
consequence we show that every vector bundle over a classifying space
$B\pi$ which can be obtained from an asymptotic representation of a
discrete group $\pi$ can be obtained also from a representation of the
group $\pi\times{\bf Z}$ into the Calkin algebra.
We give also a
generalization of the notion of Fredholm representation and show that
asymptotic representations can be viewed as asymptotic Fredholm
representations.

\end{abstract}

\section{Asymptotic representations as representations into the Calkin
algebra}\label{sec1}

Let $\pi$ be a discrete finitely presented group, and let
$F\subset\pi$ be a finite subset. Denote by
$U(n)$ the unitary group of dimension $n$ and fix a number
$\varepsilon>0$.

\begin{dfn}
{\rm
A map $\sigma:\pi\longrightarrow U(n)$ is called an {\em
$\varepsilon$-almost
representation} with respect to $F$
if $\sigma(g^{-1})=\sigma(g)^{-1}$ holds for all $g\in\pi$
and if
 $$
\norm{\sigma}_F=\sup\{\norm{\sigma(gh)-\sigma(g)\sigma(h)}:g,h,gh\in
F\}\leq\varepsilon.
 $$
}
\end{dfn}

Let $\{n_k\}$ be a strictly increasing sequence of positive integers
and let
$\sigma=\{\sigma_k:\pi\longrightarrow U(n_k)\}$ be a sequence of
$\varepsilon_k$-almost representations.
We assume that the groups $U(n_k)$ are embedded into the groups $U(n_{k+1})$
in the standard way, so it makes possible to compare almost
representations for different $k$. Then we can consider the maps
$\sigma_k\oplus 1:\pi\longrightarrow U(n_k)\oplus
U(n_{k+1}-n_k)\longrightarrow U(n_{k+1})$, which we also denote by
$\sigma_k$.

\begin{dfn}
{\rm
A sequence of $\varepsilon_k$-almost representations is called an
{\em asymptotic representation} of the group $\pi$ $($with respect to the
finite subset $F$ and a sequence $\{n_k\}$$)$ if the sequences
$\varepsilon_k$ and $\norm{\sigma_k(g)-\sigma_{k+1}(g)}:g\in F\subset\pi$
tend to zero.
}
\end{dfn}

It was shown in \cite{mis-noor} that this definition is equivalent
to the definition given in \cite{c-g-m,c-hig}. It does not depend on
the choice of the finite subset $F$ when $F$ is big enough and sets
of generators and relations of the group $\pi$ are finite.

\smallskip
Two asymptotic representations $\sigma_0$ and $\sigma_1$ are called
homotopic if there exists a family of asymptotic representations
$\sigma_t=\{\sigma_{t,k}\}$ such that
the functions $\sigma_{t,k}(g)$ are continuous for all $g\in\pi$ and
$\lim_{k\to\infty}\max_{t}\norm{\sigma_{t,k}}_F=0$.
It can be easily seen that every asymptotic representation is homotopic to
some asymptotic representation corresponding to a given sequence $\{n_k\}$.

\medskip
It is well known that the asymptotic representations are exact
representations in some more compound $C^*$-algebras~\cite{c-hig,loring-p}.
Remind this construction. Let $M_n$ be the $n\times n$ matrix algebra.
Consider the $C^*$-algebra $B=\prod_{k=1}^\infty M_{n_k}$ of norm-bounded
sequences of matrices. We suppose that the sequence $n_k$ is strictly
increasing. Denote by $B^+$ the $C^*$-algebra $B$ with adjoined unit.
Both algebras $B$ and $B^+$ contain a $C^*$-ideal $I=\oplus_{k=1}^\infty
M_{n_k}$ of sequences of matrices with norms tending to zero.

Denote the corresponding quotient algebras by $Q=B/I$ and $Q^+=B^+/I$.
Let $\bar\alpha:B^{(+)}\longrightarrow B^{(+)}$ be the right shift,
$\bar\alpha(m_1,m_2,\ldots)=(0,m_1,m_2,\ldots)$, $(m_i)\in B^{(+)}$.
As $\bar\alpha(I)\subset I$, so the homomorphism $\bar\alpha$
induces the homomorphism $\alpha:Q^{(+)}\longrightarrow Q^{(+)}$.
Let
$$
Q^{(+)}_\alpha=\{q\in
Q^{(+)}: \alpha(q)=q\}\subset Q^{(+)}
$$
be the $\alpha$-invariant
$C^*$-subalgebra. The adjoined unit gives the splittable short exact
sequence
\begin{equation}\label{exact}
Q_\alpha\longrightarrow Q^+_\alpha\longrightarrow{\bf C}.
\end{equation}
Let $e=(e_k)\in B$ be the sequence of diagonal matrices having
unities at the first place and zeroes at the other places,
$e_k\in M_{n_k}$.

\begin{lem}\label{lem}
The group $K_0(Q^+_\alpha)$ is isomorphic to ${\bf Z}\oplus{\bf Z}$
with generators $[e]$ and $[1]$;

$K_0(Q_\alpha)\cong {\bf Z}$ with generator $[e]$;

$K_1(Q_\alpha)=K_1(Q^+_\alpha)=0$.

\end{lem}

{\bf Proof.}
Let $p\in M_r(Q_\alpha)$ be a projection and let $p'\in M_r(B)$ be a
projection which is the lift of $p$. Then $p'=(p_k)$ is a sequence of
projections. This sequence is $\alpha$-invariant, hence the rank of
$p_k$ is constant beginning from some $k$. Therefore
the projection $p$ is a multiple of the projection $e$ for big enough
$k$. The case of the group $K_1(Q_\alpha)$ can be considered by the similar
way.$\quad\bullet$

\medskip
Denote by ${\cal Q}$ the Calkin algebra.

\begin{thm}\label{hom}
There exists a monomorphism
 \begin{equation}\label{psi}
\psi:Q^+_\alpha\otimes C(S^1)\longrightarrow{\cal Q}
 \end{equation}
such that the induced homomorphism
 \begin{equation}\label{psi*}
\psi_*:K_1(Q_\alpha\otimes C(S^1))\longrightarrow K_1({\cal Q})
 \end{equation}
is an isomorphism.

\end{thm}

{\bf Proof.} We start with defining the homomorphism
(\ref{psi}).
Let $V_k$ be a $n_k$-dimensional Hilbert space where the algebra
$M_{n_k}$ acts, with a fixed basis $\{e_k^{(i)}\}$, $i\leq n_k$.
Fix the embeddings
 $$
j_k:V_{k-1}\longrightarrow V_{k},
 $$
mapping the vectors
$e_{k-1}^{(i)}$ into $e_k^{(i)}$. Denote by $W_{k-1}$ the orthogonal
complements,
 $$
V_k=j_k(V_{k-1})\oplus W_{k-1}.
 $$
For $k=1$ let $W_0=V_1$. Put
 $$
H_k=V_k\oplus(\oplus_m W_{k,m});\quad W_{k,m}\cong W_k,\ m\in {\bf N};
\quad H=\oplus H_k
 $$
and define a homomorphism $\psi$ of the $C^*$-algebra $Q^+_\alpha$ into
the Calkin algebra ${\cal Q}$ of the Hilbert space $H$. Let $q\in
Q^+_\alpha$ and let $q'=(q_k)\in B^+$ be a lift of $q$.
Denote by $(q_k)_{ij}$ the matrix elements of the matrix $q_k$.
By definition the limits $\lim_{i\to\infty}(q_k)_{ii}$ exist and are
equal to each other for all $k$. Denote the value of these limits by
$\lambda(q)$. Define $\psi(q)$ as an operator on
$H$ which acts on the spaces $V_k$ by multiplication by
the matrices $q_k$, and on the spaces $W_{k,m}$ by multiplication by
$\lambda(q)$. Such operator is defined up to compact operators, hence it
gives an element $\psi(q)\in{\cal Q}$.
Consider now the $C^*$-algebra $Q^+_\alpha\otimes C(S^1)$ of continuous
$Q^+_\alpha$-valued functions over a circle.
The construction described above allows us to define a homomorphism
$Q^+_\alpha\otimes C(S^1)\longrightarrow {\cal Q}$.
Put $\psi(q\otimes 1_{C(S^1)})=\psi(q)$. It remains now to
define the image of the function $u=e^{2\pi i t}$
from $C(S^1)$. Define a Fredholm operator $F$ (with zero index) on the
Hilbert space $H$ by the following conditions:
 \begin{enumerate}
\item
$F$ maps every subspace $V_k$ (except $V_1$) onto the subspace
$V_{k-1}\oplus W_{k-1,1}$ so that $F(e_k^{(i)})=e_{k-1}^{(i)}$ for
$i<n_{k-1}$;
\item
$F$ isomorphically maps the subspace $V_1$ onto $W_{0,1}$;
\item
$F$ isomorphically maps every subspace $W_{k,m}$ onto the subspace
$W_{k,m+1}$.
\end{enumerate}
Thus defined operator $F$ commutes with the image $\psi(Q^+_\alpha)$
modulo compacts, hence we get the needed homomorphism $\psi$ if we put
$\psi(1_{Q^+_\alpha}\otimes e^{2\pi i t})=F$.

Denote the basis vectors of the subspaces $W_{k,m}$ by $w_{k,m}^{(j)}$
with $n_{k-1}+1\leq j\leq n_k$. Then we can represent the action of the
operator $F$ by the diagram

$$
\begin{array}{cccccccccccccccc}
\ldots&\lar&w_{1,2}^{(1)}&\lar&w_{1,1}^{(1)}&\lar&|&e_1^{(1)}&\lar&&e_2^{(1)}&
\lar&&e_3^{(1)}&\lar&\ldots\\
\ldots&\ldots&\ldots&\ldots&\ldots&\ldots&|&\ldots&\ldots&&\ldots&\ldots&&
\ldots&\ldots&\ldots\\
\ldots&\lar&w_{1,2}^{(n_1)}&\lar&w_{1,1}^{(n_1)}&\lar&|&e_1^{(n_1)}&
\lar&&e_2^{(n_1)}&\lar&&e_3^{(n_1)}&\lar&\ldots\\
&&&&&&&\mbox{---}\quad\mbox{---}&\mbox{---}&&&&&&&\\
\ldots&\lar&w_{2,3}^{(n_1+1)}&\lar&w_{2,2}^{(n_1+1)}&\lar&&w_{2,1}^{(n_1+1)}&
\lar&|&e_2^{(n_1+1)}&\lar&&e_3^{(n_1+1)}&\lar&\ldots\\
\ldots&\ldots&\ldots&\ldots&\ldots&\ldots&&\ldots&\ldots&|&\ldots&\ldots&&
\ldots&\ldots&\ldots\\
\ldots&\lar&w_{2,3}^{(n_2)}&\lar&w_{2,2}^{(n_2)}&\lar&&w_{2,1}^{(n_2)}&
\lar&|&e_2^{(n_2)}&\lar&&e_3^{(n_2)}&\lar&\ldots\\
&&&&&&&&&&\mbox{---}\quad\mbox{---}&\mbox{---}&&&&\\
\ldots&\lar&w_{3,4}^{(n_2+1)}&\lar&w_{3,3}^{(n_2+1)}&\lar&&w_{3,2}^{(n_2+1)}&
\lar&&w_{3,1}^{(n_2+1)}&\lar&|&e_3^{(n_2+1)}&\lar&\ldots\\
\ldots&\ldots&\ldots&\ldots&\ldots&\ldots&&\ldots&\ldots&&\ldots&\ldots&|&
\ldots&\ldots&\ldots\\
\ldots&\lar&w_{3,4}^{(n_3)}&\lar&w_{3,3}^{(n_3)}&\lar&&w_{3,2}^{(n_3)}&
\lar&&w_{3,1}^{(n_3)}&\lar&|&e_3^{(n_3)}&\lar&\ldots\\
\ldots&\ldots&\ldots&\ldots&\ldots&\ldots&&\ldots&\ldots&&\ldots&\ldots&&
\ldots&\ldots&\ldots
\end{array}
$$

It remains to show that the homomorphism (\ref{psi*}) is an isomorphism.
As $K_1(Q^+_\alpha)=0$, so by the K\"unneth formula we have
 $$
K_1(Q^+_\alpha\otimes C(S^1))\cong K_0(Q^+_\alpha)\otimes K^1(S^1).
 $$
Let $[u]\in K^1(S^1)$ be a generator. Then the group
$K_1(Q^+_\alpha\otimes C(S^1))$ is generated by two elements:
$[e]\otimes [u]$ and $[1]\otimes [u]$.
As
 $$
\psi(1\otimes u)=F\quad {\rm and}\quad {\rm ind}\, F=0,
 $$
so we obtain
$\psi_*([1]\otimes [u])=0$. Compute $\psi_*([e]\otimes [u])$. By
 \cite{connes} one has
 $$
[e]\otimes [u]=[(1-e)\otimes 1+e\otimes u],
 $$
hence
 \begin{eqnarray*}
\psi_*([e]\otimes [u])&=&\psi_*([e]\otimes [u])
=\psi_*([(1-e)\otimes 1+e\otimes
u])\\&=&[1-\psi(e\otimes 1)+\psi(e\otimes u)]=
[1-\psi(e\otimes 1)+\psi(e\otimes 1)F].
 \end{eqnarray*}
The action of the operator
 $$
F'=1-\psi(e\otimes 1)+\psi(e\otimes 1)F
 $$
is given by the formula
 $$
F'(e_k^{(j)})=\left\lbrace
  \begin{array}{ll}
e_k^{(j)}& {\rm for}\ \ j>1,\\
e_{k-1}^{(j)}& {\rm for}\ \ j=1,k>1,\\
0 &{\rm for}\ \ j=k=1;
  \end{array}\right.
\qquad
F'(w_{k,m}^{(j)})=w_{k,m}^{(j)},
$$
so the image of $F'$
coincides with the whole space $H$,
and the kernel of $F'$ is one-dimensional and is generated by the
vector $e_1^{(1)}$. Therefore
${\rm ind}\, F'=1$, hence the homomorphism
$\psi_*$ maps the generator of the group
$K_1(Q_\alpha\otimes C(S^1))\cong{\bf Z}$ into the generator of the group
$K_1({\cal Q})$.$\quad\bullet$

\medskip
Denote by ${\cal R}_a(\pi)$ the Grothendieck group of virtual asymptotic
representations of the group $\pi$. Let $\til{\cal R}_a(\pi)$ denote the
kernel of the map
 $$
{\cal R}_a(\pi)\longrightarrow {\cal R}_a(e)\cong {\bf Z},
 $$
defined by the trivial representation (here $e$ is the trivial group).
Let $B\pi$ be the classifying space for the group
$\pi$. Remember that a construction was defined
in \cite{mis-noor}, which allows to obtain a vector bundle over $B\pi$
starting from an asymptotic representation. This construction gives a
homomorphism
 \begin{equation}\label{phi}
\phi:\widetilde{\cal R}_a(\pi)\longrightarrow K^0(B\pi)
 \end{equation}
which can be described as follows.
Let $C^*[\pi]$ be the group $C^*$-algebra of the group $\pi$ and let $\xi\in
K^0_{C^*[\pi]}(B\pi)$ be the universal bundle. An asymptotic representation
$\sigma$ defines a homomorphism
 $$
\ov{\sigma}:C^*[\pi]\longrightarrow Q^+_\alpha,
 $$
which maps the universal bundle $\xi$ into some element
$\ov{\sigma}_*(\xi)\in K^0_{Q^+_\alpha}(B\pi)$.
This homomorphism defines the homomorphism
 $$
\phi':{\cal R}_a(\pi)\longrightarrow K^0_{Q^+_\alpha}(B\pi).
 $$
It follows from  the lemma \ref{lem}, from the K\"unneth formula and from
(\ref{exact}) that the lower line of the diagram
 $$
\begin{array}{ccccc}
\til{\cal R}_a(\pi)&\longrightarrow&{\cal R}_a(\pi)&
\longrightarrow&{\cal R}_a(e)\\
\downarrow&&\downarrow\lefteqn{\phi'}&&\downarrow\lefteqn{\phi'_e}\\
K^0_{Q_\alpha}(B\pi)&\longrightarrow&K^0_{Q^+_\alpha}(B\pi)&
\longrightarrow&K^0_{\bf C}(B\pi)
\end{array}
 $$
is exact, therefore the left vertical arrow is well-defined.
As we have
 $$
K^0_{Q_\alpha}(B\pi)=K^0(B\pi)\otimes K_0(Q_\alpha)\cong K^0(B\pi),
 $$
so we can define $\phi$ to be this left vertical arrow after identifying
$K^0_{Q_\alpha}(B\pi)$ with $K^0(B\pi)$.
Notice that the image of the
homomorphism $\phi'_e$ coincides with the subgroup in
$K^0(B\pi)$ generated by the trivial representations.

\smallskip
As there exists a natural isomorphism
 $$
j:K^0_A(B\pi\times S^1\times S^1)\widetilde{\longrightarrow}
K^0_{A\otimes C(S^1)}(B\pi\times S^1)
 $$
for any $C^*$-algebra $A$, so multiplication by the Bott generator
$\beta\in K^0(S^1\times S^1)$ defines an inclusion
 \begin{equation}\label{..}
\begin{array}{ccccl}
\ov{\beta}:K^0_A(B\pi)&\stackrel{\otimes\beta}{\longrightarrow}&
K^0_A(B\pi\times S^1\times S^1)&\stackrel{j}{\longrightarrow}&
K^0_{A\otimes C(S^1)}(B\pi\times S^1)\\
&&&=&K^0_{A\otimes C(S^1)}
(B(\pi\times{\bf Z})).
\end{array}
 \end{equation}
In the case $A=Q_\alpha$ we will write $\ov{\beta}$ instead of
$\ov{\beta}_{Q_\alpha}$.

Now denote by ${\cal R}_{\cal Q}(\pi)$ the group of (virtual)
representations of the group $\pi$ into the Calkin algebra.
It is easily seen that the homomorphism (\ref{phi}) allows us to define
a homomorphism
 \begin{equation}\label{...}
\ov{\psi}:{\cal R}_a(\pi)\longrightarrow {\cal R}_{\cal Q}(\pi\times
{\bf Z})
 \end{equation}
given by the formula $\ov{\psi}(\ov{\sigma})=\psi(\ov{\sigma}\otimes id)$
for $\ov{\sigma}\in{\cal R}_a(\pi)$ and
 $$
\ov{\sigma}\otimes id: C^*[\pi\times {\bf Z}]\cong
C^*[\pi]\otimes C(S^1)\longrightarrow
Q^+_\alpha\otimes C(S^1).
 $$

Let $\eta\in K^0_{C^*[\pi\times{\bf Z}]}(B(\pi\times Z))$
be the universal bundle over $B(\pi\times{\bf Z}$.
Then there exists also a homomorphism
 \begin{equation}\label{.}
f:{\cal R}_{\cal Q}(\pi\times{\bf Z})\longrightarrow
K^0_{\cal Q}(B(\pi\times{\bf Z}))
 \end{equation}
defined as the image of $\eta$ over
$B(\pi\times{\bf Z})$ under the representations into the Calkin algebra.

All these homomorphisms (\ref{psi*}) --  (\ref{...})
can be represented by the diagram
 \begin{equation}\label{diagr}
\begin{array}{ccccc}
\widetilde{\cal R}_a(\pi)&\stackrel{\phi}{\longrightarrow}&
K^0_{Q_\alpha}(B\pi)&\stackrel{\ov{\beta}}{\longrightarrow}&
K^0_{Q_\alpha\otimes C(S^1)}(B\pi\times S^1)\\
\downarrow&&&&\downarrow\lefteqn{\psi_*}\\
{\cal R}_a(\pi)&\stackrel{\ov{\psi}}{\longrightarrow}&
{\cal R}_{\cal Q}(\pi\times{\bf Z})&\stackrel{f}{\longrightarrow}&
K^0_{\cal Q}(B\pi\times S^1),
\end{array}
 \end{equation}
where the left vertical arrow is the inclusion.

\begin{thm}\label{diagram}
The diagram $($\ref{diagr}$)$ is commutative.

\end{thm}

{\bf Proof.} Direct calculation.$\quad\bullet$

\medskip
So now we have the homomorphism
 \begin{equation}\label{*}
\ov{\phi}:\widetilde{\cal R}_a(\pi)\longrightarrow
K^0_{\cal Q}(B\pi\times S^1)\cong K^1(B\pi\times S^1).
 \end{equation}

\smallskip
Theorem \ref{diagram} shows that every element of the group $K^0(B\pi)$
which can be obtained by an
asymptotic representation of the fundamental group $\pi$ can be obtained
also by a
representation of the group $\pi\times {\bf Z}$ into the Calkin algebra.

\section{Asymptotic representations as Fredholm
representations}\label{sec2}
\setcounter{equation}{0}

Notice that the group $K^0_{\cal Q}(B\pi\times S^1)$ can be decomposed
into direct sum:
 \begin{equation}\label{sum}
K^0_{\cal Q}(B\pi\times S^1)=K^0_{\cal Q}(B\pi)\oplus
K^0_{\cal Q}(B\pi\wedge S^1)\cong K^1(B\pi)\oplus K^0(B\pi)
 \end{equation}
induced by an inclusion map $i:s_0\longrightarrow S^1$, where $s_0\in S^1$,
and the image of the homomorphism $\ov{\phi}$ (\ref{*}) lies only
in the second summand of (\ref{sum}). Indeed, consider the composition
of the map $\ov{\phi}$ with the map $i^*:K^0_{\cal
Q}(B\pi\times S^1)\longrightarrow K^0_{\cal Q}(B\pi)$. But as the
multiplication by the Bott generator is involved in the map $\ov{\phi}$,
so its composition with $i^*$ gives the zero map. Therefore the image of
the group $\widetilde{\cal R}_a(\pi)$ lies in $K^0(B\pi)$ and hence
defines a (virtual) vector bundle over $B\pi$.

\smallskip
On the other hand the image of the map (\ref{.}) need not be
contained in the second summand of (\ref{sum}). It would be so if the
representation of the group $\pi\times {\bf Z}$ into the Calkin algebra
would be a part of a Fredholm representation of the group $\pi$
\cite{mis}. Using the notion of asymptotic representations we can now give
a generalization of the Fredholm representations which would also ensure
that the image of such representations would lie in the second summand of
(\ref {sum}).

\medskip
Let $\rho:\pi\times {\bf Z}\longrightarrow {\cal Q}$ be a representation
into the Calkin algebra and let $F\subset\pi$ be a finite subset of $\pi$
containing its generators. Denote by $B(H)$ the algebra of bounded
operators on a separable Hilbert space $H$. Let $q:B(H)\longrightarrow
{\cal Q}$ be the canonical projection.

 \begin{dfn}
{\rm
We call a map
$\tau:\pi\longrightarrow B(H)$ an {\em $\e$-trivialization} for $\rho$ if
\begin{enumerate}
\item
$\|\tau(gh)-\tau(g)\tau(h)\|\leq\e$ for any $g,h\in F\subset\pi$,
\item
$q(\tau(g))=\rho(g;0)$ for any $g\in\pi$, $(g;0)\in\pi\times{\bf Z}$.
\end{enumerate}
}
 \end{dfn}

 \begin{dfn}
{\rm
Suppose that for every $\e>0$ there exists an $\e$-trivialization
$\tau_\e$ for $\rho$. Then the pair $(\tau_\e,\rho)$ is called an
{\em asymptotic Fredholm representation}.
}
 \end{dfn}

Let $u$ be a generator of the group $\bf Z$.
Notice that the image of the group $\pi$ under $\varepsilon$-trivializations
commutes with
some Fredholm operator $F=\rho(0,u)$ modulo
compacts.
Denote the group of all asymptotic Fredholm representations by ${\cal
R}_{aF}(\pi)\subset {\cal R}_{\cal Q}(\pi\times{\bf Z})$.

 \begin{prop}
The image of ${\cal R}_{aF}$ under the map $f$ $($\ref{.}$)$
lies in the group $K^0(B\pi\wedge
S^1)$.
 \end{prop}

{\bf Proof.}
It was described in the paper \cite{mis-noor} how to construct a bundle
over $B\pi$ with the fibers isomorphic to the Hilbert space $H$ and with
the structural group $GL(H)$
starting from an almost representation $\tau_\e$ for small
enough $\e$.
To do so one should construct transition functions acting on the fibers
$\xi_x$,
 $$
T_g(x):\xi_x\longrightarrow\xi_{gx}
 $$
for $g\in\pi$, $x\in E\pi$.
One should chose representatives $\{a_\a\}$ in each orbit of the set of
vertices of $E\pi$ and define $T_g(a_\a)=\tau_\e(g)$. Take now an
arbitrary vertex $b\in E\pi$. Then there exists such $h\in \pi$ that
$b=h(a_\a)$ and we should put $T_g(b)=\tau_\e(gh)\tau_\e^{-1}(h)$.
Further these transition functions should be extended by linearity to all
simplexes of $E\pi$.
But obviously $q(T_g(b))=q(\tau_\e(gh)\tau_\e^{-1}(h))=q(\tau_\e(g))$,
hence after we pass to quotients,
the transition functions $q(T_g(x))$
would become constant and the bundle with the
structural group being the invertibles of the Calkin algebra. But as any
bundle with fibers $H$ is trivial, so the quotient bundle with fibers
isomorphic to the Calkin algebra is trivial too and the projection of
$\varphi({\cal R}_{aF})$ onto the first summand of $K^0_{\cal Q}(B\pi)$ is
equal to zero.$\quad\bullet$

\medskip
{\bf Remark.}
In the section \ref{sec1} we have seen that an asymptotic
representation $\sigma=(\sigma_n)$ defines a homomorphism
$\rho:\pi\times{\bf Z}\longrightarrow GL({\cal Q})$ into the group of
invertibles of the Calkin algebra. But the same asymptotic representation
gives $\e$-trivializations of $\rho$ for any $\e>0$. Indeed we can put
 $$
\tau_{\e_n}(g)=(1,1,\ldots,1,\sigma_n(g),\sigma_{n+1}(g),\ldots).
 $$
So asymptotic representations define asymptotic Fredholm representations
and we have an inclusion ${\cal R}_a(\pi)\subset{\cal R}_{aF}(\pi)$.


\section{Asymptotic representations and extensions}\label{sec3}
\setcounter{equation}{0}

It was shown in \cite{c-hig} that using a quasi central approximate unity
one can construct an aymptotic representation out of an extension of
$C^*$-algebras. We study how this construction is related to asymptotic
representations of the initial $C^*$-algebra.

\smallskip
Let $A$ be a $C^*$-algebra such that there exists a homomorphism
$q_A:A\longrightarrow {\bf C}$ of $A$ into the complex numbers and let
$F\subset A$ be a finite set of generators for $A$. We can repeat our
definitions of asymptotic representations from the section \ref{sec1} in
the case of $C^*$-algebras instead of discrete groups. For smplicity sake
we assume that dimension of an almost representation $\sigma_n$ equals to
$n$.

\smallskip
We consider a sequence of maps $\til{\sigma}_n:A\longrightarrow M_n$
where $M_n$ acts on a fi\-ni\-te-di\-men\-si\-onal Hilbert space $V_n$.
Take an in\-fi\-ni\-te-di\-men\-si\-onal Hilbert space $H_n\supset V_n$
and define a map $\sigma_n:A\longrightarrow B(H)$ by $\sigma_n(a)=
\til{\sigma}_n(a)\oplus q_A(a)$ where $a\in A$ and $q_A(a)$ is a scalar on
$V_n^\perp$. The map $\sigma=\oplus_n\sigma_n$ gives an asymptotic
representation of the $C^*$-algebra $A$ if (after identifying all $H_n$)
the norms
$\|\sigma_{n+1}(a)-\sigma_n(a)\|$ tend to zero for $a\in F$.

\smallskip
Consider the $C^*$-algebra ${\cal E}$ in the Hilbert space $\oplus_n H_n$
generated by $\sigma(a)$, $a\in A$ and by the translation operator $F$
defined in the proof of the theorem \ref{hom}. Then one has a short exact
sequence
 \begin{equation}\label{exactseq}
\oplus_n M_n\longrightarrow {\cal E}\longrightarrow A\otimes C(S^1).
 \end{equation}
We consider a discrete version of the construction of \cite{c-hig}.
If $e_n\in\oplus_n M_n$ is a quasi central approximate unity \cite{arv}
then the exact sequence (\ref{exactseq}) defines an asymptotic representation
 \begin{equation}\label{as-connes}
\rho_n:A\otimes C(S^1)\otimes C_0(S^1)
\stackrel{as}{\longrightarrow} \oplus_n M_n
 \end{equation}
given by the formula
 $$
\rho_n(a\otimes g\otimes f)=(a\otimes g)'\cdot f(e_n),
 $$
where $(a\otimes g)'\in {\cal E}$ is a lift for $a\otimes g$.
Denote by $K_n\subset\oplus_nH_n$ the subspace on which one has $e_n\neq 0$,
$e_n\neq 1$.
If $K_n$ is
fi\-ni\-te-di\-men\-si\-onal, $k(n)=\dim K_n$,
then we can get a sequence of
fi\-ni\-te-di\-men\-si\-onal almost representations
 $$
\ov{\rho}_n:A\otimes C(S^1)\otimes C_0(S^1)
\stackrel{as}{\longrightarrow} \oplus_n M_{k(n)}
 $$
by the formula
 $$
\ov{\rho}_n(a\otimes g\otimes f)=P_n(a\otimes g)'P_n\cdot f(e_n),
 $$
where $P_n$ is the projection onto $K_n$, $g\in C(S^1)$, $f\in C_0(S^1)$,
$f(0)=0$.

\medskip
There exists also a well-known asymptotic representation $\beta=(\beta_m)$ of
the commutative $C^*$-algebra $C(S^1\times S^1)$ into the matrix algebra
$M_n$ given by
 $$
\beta_m(e^{2\pi ikx}e^{2\pi ily})=T_m^kU_m^l, \qquad
\beta_m:C(S^1)\otimes C(S^1)\stackrel{as}{\longrightarrow} M_m,
 $$
where $T_m$ and $U_m$ are the $m$-dimensional
Voiculescu matrices, $T_m$ is a translation and
$U_m$ is a diagonal matrix \cite{voi}.
This asymptotic representation
realizes the Bott isomorphism.

\smallskip
The tensor product of $\til{\sigma}$ by $\beta$ gives an asymptotic
representation
 $$\
\til{\sigma}_n\otimes\beta_m:A\otimes C(S^1)\otimes C_0(S^1)
\stackrel{as}{\longrightarrow}M_{n\times m}.
 $$
Put $\widehat{\sigma}=P_n\sigma_n P_n$

 \begin{thm}
There exists a quasi central approximate unity $e_n$ for the exact
sequence $($\ref{exactseq}$)$
and a sequence of numbers $m(n)$
such that the asymptotic representations
$\rho_n$ $($\ref{as-connes}$)$ and
$\widehat{\sigma}_n\otimes\beta_{m(n)}$
into the matrix algebra $M_{(n+m(n))\times m(n)}$ are
equivalent. Moreover for $a\in F\subset A$, $e^{2\pi ix}\in C(S^1)$ and
$f\in C_0(S^1)$ the norms
 $$
\|\rho_n(a\otimes e^{2\pi ix}\otimes f)-
\widehat{\sigma}_n(a)\otimes\beta_{m(n)}(e^{2\pi ix}\otimes f)\|
 $$
tend to zero.
 \end{thm}

{\bf Proof.}
Define a sequence of integers $m(n)$ so that the following conditions
would be satisfied:
\begin{enumerate}
\item
$\max_{n\leq k\leq n+m(n)} \|\sigma_n(a)-\sigma_{n+k}(a)\|$ tend to zero for
any $a\in F\subset A$,
\item
$m(n)$ tends to infinity.
\end{enumerate}

As we now assume that ${\rm dim}\, V_n=n$, so ${\rm dim}\, W_{k,l}=1$.
We denote the basis vectors of $V_n$ by $v_n^{(j)}$, $j=1,\ldots, n$ and
the vectors of the subspaces $W_{k,l}$ by $v_{k-l}^{(k)}$.
In these denotations
the shift operator $F$ constructed in the section \ref{sec1} acts by
the simple formula $F v_n^{(j)}=v_{n-1}^{(j)}$.

\smallskip
Put
 $$
a_{n,i}^{(j)}=\left\lbrace
\begin{array}{ll}
1,&{\rm for}\ \  -n\leq i\leq n,\\
\frac{m(n)-i+n}{m(n)},&{\rm for}\ \  n<i\leq n+m(n) \ \ {\rm and}\ \ j\leq
n+m(n),\\
0, & {\rm otherwise}
\end{array}\right.
 $$
and define the diagonal quasi central approximate unity $e_n$ by the
formula
 $$
e_n v_i^{(j)}=a_{n,i}^{(j)}v_i^{(j)}.
 $$
Notice that $e_n$ exactly commutes with $\sigma_n(a)$ (which are diagonal
too) and for any $f\in
C_0(S^1)$, $f(0)=0$, $f(e_n)$ almost commutes with $F$.

Obviously the norms
 $$
\| \rho_n(a\otimes e^{2\pi ix}\otimes f)-
\ov{\rho}_n(a\otimes e^{2\pi ix}\otimes f)\|
 $$
tend to zero, so we may deal with finite-dimensional asymptotic
representations $\ov{\rho}_n$ instead of $\rho_n$.
Direct calculations shows that
 $$
\ov{\rho}_n(a\otimes e^{2\pi ix}\otimes f)=P_na'FP_n f(e_n)=
P_n(\oplus_{i=n+1}^{m(n)}\sigma_i)(a)P_nFP_n f(e_n)=
\oplus_{i=n+1}^{m(n)}\widehat{\sigma}_i(a)FP_nf(e_n),
 $$
where $a'=\oplus_{i=1}^\infty\sigma_i(a)\in{\cal E}$ (resp. $F$) is a lift
for $a$ (resp. for $e^{2\pi ix}$),
and $f(e_n)$ is the diagonal matrix with elements $\frac{m(n)-i+n}{m(n)}$
on the diagonal. It is easily seen that the norms
 $$
\|P_n F P_n f(e_n)-
T_{m(n)} f(e_n)\|
 $$
tend to zero and that
 $$
T_{m(n)}
f(e_n)=\beta_{m(n)}(e^{2\pi ix}\otimes f).
 $$
So it remains to notice that
by our choice of the sequence $m(n)$ the norms
 $$
\|
(\oplus_{i=n+1}^{m(n)}\widehat{\sigma}_i)(a)-
(\oplus_{i=n+1}^{m(n)}\widehat{\sigma}_n)(a)\|
 $$
tend to zero too.$\quad\bullet$

\medskip
{\bf Remark.}
Unfortunately in the case of discrete group $C^*$-algebras we do not know
if it is possible to construct a natural map from the group
${\cal R}_a(C^*(\pi)\otimes
C(S^1)\otimes C_0(S^1))$ into $K^0(B\pi)$ which would extend the map
(\ref{phi}) and close the diagram
 $$
\begin{array}{rcccl}
\widetilde{\cal R}_a(\pi)\!\!\!\!\!\!\!\!\!\!&&\longrightarrow&&
\!\!\!\!\!\!\!\!\!\!\!\!{\cal R}_a(C^*[\pi]
\otimes C(S^1)\otimes C_0(S^1))\\
&{\scriptstyle \phi}\!\!\!\searrow&&\swarrow\!\!\!
{\scriptstyle ?}&\\
&&K^0(B\pi).&&
\end{array}
 $$

\vspace{2.5cm}
\noindent
\parbox{7cm}{V. M. Manuilov\\
Dept. of Mechanics and Mathematics\\
Moscow State University\\
Moscow, 119899, RUSSIA\\
e-mail: {\tt manuilov@mech.math.msu.su}}
\hfill
\parbox{7cm}{A. S. Mishchenko\\
Dept. of Mechanics and Mathematics\\
Moscow State University\\
Moscow, 119899, RUSSIA\\
e-mail: {\tt asmish@mech.math.msu.su}}

\end{document}